# REVISITING GUO–HUANG–WANG–WEI POINTWISE ESTIMATES FOR RADIAL SOLUTIONS TO THE FOURTH-ORDER HÉNON EQUATION IN THE PUNCTURED SPACE

QUỐC ANH NGÔ AND TRUNG NGUYEN

ABSTRACT. We study pointwise estimates for positive, $C^4$-solution $u$ to the fourth-order Hénon equation

$$\Delta^2 u = |x|^\sigma u^p \quad \text{in } \mathbf{R}^n \backslash \{0\}$$

with $n \geqslant 5$, $p > 1$, and $\sigma > -4$. In an earlier work published in *J. Anal. Math.* **140** (2020) 371–394, Z. Guo, X. Huang, L. Wang, and J. Wei established the following pointwise estimate

$$-\Delta u \geqslant \frac{2}{n-4}\frac{|\nabla u|^2}{u} + \sqrt{\frac{n(n-4)}{n^2-4}}\, u^{\frac{n}{n-4}} \quad \text{in } \mathbf{R}^n \backslash \{0\}$$

for radial solutions to the critical, fourth-order Lane–Emden equation, namely $\sigma = 0$ and $p = (n+4)/(n-4)$. Inspired by this interesting result, we establish a similar pointwise inequality

$$-\Delta u \geqslant \frac{2}{n-4}\frac{|\nabla u|^2}{u} + \sqrt{\frac{2}{p+1-\frac{8}{n(n-4)}}}\, |x|^{\frac{\sigma}{2}} u^{\frac{p+1}{2}} \quad \text{in } \mathbf{R}^n \backslash \{0\}$$

for radial solutions to the fourth-order Hénon equation in $\mathbf{R}^n \backslash \{0\}$ for any $\sigma \geqslant 0$ and $p \geqslant (n+4+\sigma)/(n-4)$. We also show that the constant $2/(n-4)$ is sharp and that the thresholds $\sigma \geqslant 0$ and $p \geqslant (n+4+\sigma)/(n-4)$ are necessary. As in our previous study on the case of $\mathbf{R}^n$, our proof is also based on the Bernstein technique via the maximum principle, hence appears to be simpler and more constructive. However, to avoid possible singularity at the origin, we develop a Bernstein technique in one-dimensional space, which is capable of handling equations in punctured spaces.

## 1. INTRODUCTION

In this work, we continue our study on pointwise estimates for non-negative, non-trivial, $C^4$-solution $u$ to the fourth-order Hénon equation

$$\Delta^2 u = |x|^\sigma u^p \qquad (1.1)_\sigma$$

in the punctured space $\mathbf{R}^n \backslash \{0\}$ with $n \geqslant 5$, $\sigma \geqslant 0$. There are reasons why we are interested in $(1.1)_\sigma$. In fact, equations of the form $(1.1)_\sigma$ have been studied extensively since they have their roots in various branches of applied and pure mathematics. For e.g., the simplest model $(1.1)_0$ is known as the fourth-order Lane–Emden equation generalized the second-order case, namely then second-order Lane–Emden equation

$$-\Delta u = u^p,$$

*Date*: **5th Nov, 2026 at 00:22.**

which has its root in astrophysics, see [Hen73], and seems to be relevant in Yang–Mills–Higgs theory, see [Avi87]. See also [GG06], [GGS10], and [GN22]. In addition, singular solutions to the second-order Lane–Emden equation and its generalization

$$-\Delta u = |x|^\sigma u^p \quad \text{in } \mathbf{R}^n$$

with undefined sign of $\sigma$ known as the second-order Hardy–Hénon equation, have some role in the scalar curvature problem with prescribed singular set in conformal geometry; see [MP96] and the reference therein.

In our previous work [NN26b], we successfully obtain a new proof of the following interesting pointwise inequality

$$-\Delta u \geqslant \frac{2}{n-4}\frac{|\nabla u|^2}{u} + \sqrt{\frac{2}{p+1-\frac{8}{n(n-4)}}}|x|^{\frac{\sigma}{2}} u^{\frac{p+1}{2}} \quad \text{in } \mathbf{R}^n \tag{1.2}$$

obtained earlier by Fazly et al. in [FWX15] for bounded positive $C^4$-solutions $u$ to $(1.1)_\sigma$ in the whole space $\mathbf{R}^n$ and for $\sigma \geqslant 0$ but implicitly near 0 and in the supercritical regime $p > (n+4+2\sigma)/(n-4)$. To be more precise, we establish the validity of (1.2) for any $C^4$-solution not necessarily bounded in the regime $n \geqslant 5$, $\sigma \geqslant 0$, and $p \geqslant (n+4+2\sigma)/(n-4)$. Apparently, (1.2) can be considered as an improvement of the following Modica type estimate

$$-\Delta u \geqslant \sqrt{\frac{2}{p+1}}|x|^{\frac{\sigma}{2}} u^{\frac{p+1}{2}} \quad \text{in } \mathbf{R}^n, \tag{1.3}$$

see [Pha12, Lemma 2.7], see also [Mod85], and the following gradient estimate

$$-\Delta u \geqslant \frac{2}{n-4}\frac{|\nabla u|^2}{u} \quad \text{in } \mathbf{R}^n, \tag{1.4}$$

see [LX25] and section 4 below. To relax several restrictions in [FWX15] we use the Bernstein technique in the similar manner of what we did in [NNP18]. We notice that such a improvement is not trivial because the two inequalities (1.3) and (1.4) have different scaling invariance. Loosely speaking, while the latter inequality is invariant with respect to both homogeneity and scaling invariance $\lambda \mapsto \lambda^\theta u(\lambda\cdot)$ with $\theta = (4+\sigma)/(p-1)$, the former inequality is only invariant with respect to the above scaling invariance. This makes proofing (1.2), especially the sharpness of the constant $2/(n-4)$, is delicated.

Inspired by the work [FWX15], Guo et al. prove in [GHWW20] the following inequality

$$-\Delta u \geqslant \frac{2}{n-4}\frac{|\nabla u|^2}{u} + \sqrt{\frac{n(n-4)}{n^2-4}}\, u^{\frac{n}{n-4}} \quad \text{in } \mathbf{R}^n\backslash\{0\} \tag{1.5}$$

for solutions $u$, with a non-removable singularity at the origin, to $(1.1)_0$ in $\mathbf{R}^n\backslash\{0\}$ in the critical case $p = (n+4)/(n-4)$. (The sharpness of the coefficient $2/(n-4)$ in (1.5) is not mentioned in [GHWW20].) It is worth noting that if we set $\sigma = 0$ and $p = (n+4)/(n-4)$, the two inequalities (1.2) and (1.5) are identical. In addition, if the singularity at the origin is removable, then $u \in C(\mathbf{R}^n) \cap C^4(\mathbf{R}^n\backslash\{0\})$ in the extension sense. Hence, by standard regularity result, $u$ actually belongs to the class $C^4(\mathbf{R}^n)$. Thus, (1.5) actually holds true for any solution with or without singularity at the origin, thanks to [FWX15] and [NN26b]. As highlighted by the authors, the novelty of (1.5) is that it allows solutions to be unbounded. A typical example of such a solution is $C|x|^{-(n-4)/2}$ for suitable $C > 0$.

The goal of this paper is to slightly extend the above result to the fourth-order Hénon equation $(1.1)_\sigma$ in the punctured space $\mathbf{R}^n\backslash\{0\}$ with arbitrary $\sigma \geqslant 0$. That said, we prove the following.

**Theorem 1.1** (sharp Guo–Huang–Wang–Wei's inequality for radial functions). *Let $u$ be a non-negative, non-trivial, radial, $C^4$-solution to $(1.1)_\sigma$ in $\mathbf{R}^n\backslash\{0\}$ with $n \geqslant 5$, $\sigma \geqslant 0$, and $p \geqslant (n+4+\sigma)/(n-4)$. Then, $u > 0$ and enjoys the pointwise inequality*

$$-\Delta u \geqslant \frac{2}{n-4}\frac{|\nabla u|^2}{u} + \sqrt{\frac{2}{p+1-\frac{8}{n(n-4)}}}|x|^{\frac{\sigma}{2}}u^{\frac{p+1}{2}} \quad \text{in } \mathbf{R}^n\backslash\{0\}. \tag{1.6}$$

*Moreover, the constant $2/(n-4)$ is sharp in the sense that it cannot be replaced by any bigger number and the threshold $(n+4+\sigma)/(n-4)$ is optimal in the sense that it cannot be replaced by any smaller number.*

For the validity of Theorem 1.1, perhaps one has to test (1.6) by the standard singular solution

$$u_{\text{sing}}(x) = L|x|^{-\theta}, \quad \theta = \frac{4+\sigma}{p-1}$$

with

$$L = \big(\theta(\theta+2)(n-\theta-2)(n-\theta-4)\big)^{\frac{1}{p-1}}.$$

Then, (1.6) becomes

$$\theta(n-2-\theta) \geqslant \frac{2}{n-4}\theta^2 + \sqrt{\frac{2}{p+1-\frac{8}{n(n-4)}}}\sqrt{\theta(\theta+2)(n-2-\theta)(n-4-\theta)}.$$

Equivalently,

$$\frac{n-2}{n-4}\sqrt{\frac{\theta(n-4-\theta)}{(\theta+2)(n-2-\theta)}} \geqslant \sqrt{\frac{2}{p+1-\frac{8}{n(n-4)}}},$$

which is always true under the hypotheses $n \geqslant 5$, $\sigma \geqslant 0$, and $p \geqslant (n+4+\sigma)/(n-4)$.

To prove Theorem 1.1, we use the Bernstein technique in the similar manner of what we did in [**NNP18**] and in [**NN26b**]. See also [**CCR23**]. This technique, introduced in the local case by Bernstein himself in his pioneering papers [**Ber06**, **Ber10**] over a hundred years ago, remains a powerful tool for establishing derivative estimates, either in pointwise or in integral sense, for solutions of elliptic equations. Heuristically, this technique is done through the use of auxiliary and cut-off functions together with the maximum principle. Consequently, higher order derivatives of solutions can be estimated in terms of lower order ones, but very often the size of the reference domain should appear in these estimates as a price we pay. See [**CDV22**, page 1598] for the simplest example demonstrating this technique.

The Bernstein technique has now been studied and extended to many other situations. As a matter of fact, it is impossible to provide a comprehensive account of all the results related to the Bernstein technique. For the interested readers, we refer to [**LU64**, **Oba71**, **Ser71**, **Lio80**, **GS81**, **Bar91**, **BV91**, **FV10**, **CDV22**, ?], to name a few, and the references therein for several typical examples and their subsequent influence often in PDE theory and geometric analysis.

Theorem 1.1 deserves some long comments. First, it is apparent that radial symmetry is required. Although the original statement of the result of Guo et al. in [**GHWW20**] does not require such a strong property, but it is implicitly included due to a classical classifictation of Lin for solutions to the fourth-order Lame–Emden equation $(1.1)_0$, see [Lin98]. It seems to us that such a classification result is no longer true in the case of fourth-order Hénon equations, namely $\sigma > 0$. See [GGN13] for a relavent result just in the second-order case. Therefore, it is natural to carry out the same property in order to handle any $\sigma \geqslant 0$. Nevertheless, combined with the result in [NN26b] for the case of the whole space, we have the following corollary without any radial symmetry assumption.

**Corollary 1.2** (sharp Guo–Huang–Wang–Wei's inequality for functions in the critical case). *Let $u$ be a non-negative, non-trivial, $C^4$-solution to $(1.1)_0$ in $\mathbf{R}^n\backslash\{0\}$ with $n \geqslant 5$ and $p = (n+4)/(n-4)$. Then, $u > 0$ and enjoys the pointwise inequality*

$$-\Delta u \geqslant \frac{2}{n-4}\frac{|\nabla u|^2}{u} + \sqrt{\frac{n(n-4)}{n^2-4}}\, u^{\frac{n}{n-4}} \quad \text{in } \mathbf{R}^n\backslash\{0\}$$

*Moreover, the constant $2/(n-4)$ is sharp in the sense that it cannot be replaced by any bigger number.*

Unfortunately, even in the critical case but with $\sigma > 0$, we do not know whether or not Corollary 1.2 remains valid. However, the radial symmetry assumption is just a techinical assumption. We notice that using integral representation the weaker version (1.4) holds in $\mathbf{R}^n\backslash\{0\}$, namely

$$-\Delta u \geqslant \frac{2}{n-4}\frac{|\nabla u|^2}{u} \quad \text{in } \mathbf{R}^n\backslash\{0\},$$

see section 4. It is also worth noting that by a gluing method for fourth-order ODEs, there are infinitely many nonradial singular solutions to the fourth-order Lane–Emden equation $(1.1)_0$ in $\mathbf{R}^n\backslash\{0\}$ with $n \geqslant 6$ and

$$\frac{n+3}{n-5} < p < p_{\rm JL}(n-1),$$

where $p_{\rm JL}(n)$ is the Joseph–Lundgren exponent given by

$$p_{\rm JL}(n) = \begin{cases} \dfrac{(n-2)^2 - 4n + 8\sqrt{n-1}}{(n-2)(n-10)} & \text{if } n \geqslant 11, \\ +\infty & \text{if } 2 \leqslant n \leqslant 10; \end{cases}$$

see [NN26a, GWY16]. It would be interesting to see if such nonradial singular solutions actually fulfill (1.6).

Second, given $n \geqslant 5$, $\sigma > -4$, and $p > 1$, although $(1.1)_\sigma$ in $\mathbf{R}^n\backslash\{0\}$ admits solutions if

$$\text{either} \quad 1 < p < \frac{n+2+\sigma}{n-2} \quad \text{or} \quad p > \frac{n+\sigma}{n-4},$$

see [Ngo24], but by testing with the standard singular solution, it turns out that for the validity of (1.8), it is necessary to have $p > (n+\sigma)/(n-4)$. It fact, given $n \geqslant 5$, $\sigma > -4$, and $1 < p < (n+2+\sigma)/(n-2)$, the inequality (1.8) is reverse, namely any radial solution $u$ to $(1.1)_\sigma$ in $\mathbf{R}^n\backslash\{0\}$ enjoys

$$-\Delta u < \frac{2}{n-4}\frac{|\nabla u|^2}{u} \quad \text{in } \mathbf{R}^n\backslash\{0\}, \tag{1.7}$$

see [NN26c].

Third, given $n \geqslant 5$ and $\sigma > -2$, fixing any $p$ satisfying $(n+2+\sigma)/(n-4) < p < (n+4+\sigma)/(n-4)$, there is still no constant $C > 0$ such that

$$-\Delta u \geqslant \frac{2}{n-4}\frac{|\nabla u|^2}{u} + C|x|^{\frac{\sigma}{2}} u^{\frac{p+1}{2}} \quad \text{in } \mathbf{R}^n\backslash\{0\} \tag{1.8}$$

holds; see Remark 3.2. In other words, the threshold $(n+4+\sigma)/(n-4)$ is sharp for the validity of inequalities of the form (1.6), even (1.8), as claimed in Theorem 1.1. Fourth, we also note that, given $n \geqslant 5$, $\sigma \in (-2, 0)$, and $p > (n+4+2\sigma)/(n-4)$, there is still no constant $C > 0$ such that(1.8) holds. In other words, the condition $\sigma \geqslant 0$ is necessary for the validity of (1.6); see Remark 3.1.

From all the discussion above, we obtain the following table.

| | $-4 < \sigma \leqslant 2$ | $-2 < \sigma < 0$ | $\sigma \geqslant 0$ |
|---|---|---|---|
| $\frac{n+\sigma}{n-4} < p \leqslant \frac{n+2+\sigma}{n-4}$ | | | |
| $\frac{n+2+\sigma}{n-4} < p < \frac{n+4+\sigma}{n-4}$ | | Remarks 3.2, 3.1 | Remark 3.2 |
| $p \geqslant \frac{n+4+\sigma}{n-4}$ | | Remark 3.1 | Theorem 1.1 |

TABLE 1. On the thresholds for $\sigma > -4$ and $p > (n+\sigma)/(n-4)$.

Hence, we conclude that all the thresholds for $\sigma$ and $p$ in Theorem 1.1 are optimal. (It would be interesting if there is any counter-example for the failure of Theorem 1.1 corresponding to empty cells in Table 1 above.)

Finally, we go back to the weaker version (1.4) and prove the following result.

**Theorem 1.3** (sharp pointwise inequality). *Let $u$ be a non-negative, non-trivial, $C^4$-solution to the fourth-order Hardy–Hénon equation $(1.1)_\sigma$ in $\mathbf{R}^n\backslash\{0\}$ with $n \geqslant 5$, $\sigma > -4$, and $p > (n+\sigma)/(n-4)$. Then, $u$ is strictly positive and the following pointwise inequality*

$$-\Delta u \geqslant \frac{2}{n-4}\frac{|\nabla u|^2}{u} \quad \text{in } \mathbf{R}^n\backslash\{0\}$$

*holds. Moreover, the constant $2/(n-4)$ is sharp in the sense that it cannot be replaced by any bigger number.*

We prove Theorem 1.3 in section 4. While the inequality is more or less standard, as we mentioned earlier, the sharpness of the coefficient $2/(n-4)$ is dedicated. We are successfully in prove the sharpness of $2/(n-4)$, by using a similar idea as in [**NN26b**]. However, unlike the case considered in [**NN26b**], solutions to $(1.1)_\sigma$ in $\mathbf{R}^n\backslash\{0\}$ are not always bounded. Hence, we need to look for suitable radial solutions with appropriate decay rate at infinity.

Before mentioning the organization of this paper and closing this section, we should notice that our inequality (1.6) for both $\mathbf{R}^n$ and the punctured space $\mathbf{R}^n\backslash\{0\}$ in the critical case can be used to study the exact behavior of radial solutions to $(1.1)_\sigma$ in $\mathbf{R}^n$ and in $\mathbf{R}^n\backslash\{0\}$ with $n \geqslant 5$, $p > 1$, and $\sigma \geqslant 0$. (The corresponding case $-4 < \sigma < 0$ is already known, see e.g. the work [**FK19**].) This seems to be the first application of the inequality (1.6) in PDE. Details will appear elsewhere. The organization of this paper is as follows.

## Contents

## 2. Preliminaries

Throughout this section by $u$ we mean a $C^4$-solution to $(1.1)_\sigma$ in $\mathbf{R}^n \backslash \{0\}$ with $n \geqslant 5$, $p > 1$, and $\sigma > -4$. We spend this section to collect some results necessarily for our approach, which is essentially the classical Bernstein technique. We note that since we deal with pointwise estimate, proof of many results in this section can be found in [NN26b]. Also, these results hold true for any solution not necessarily radial. In some cases, it is required $\sigma \geqslant 0$.

That said, the following change of variables

$$k = \frac{p-1}{2} > 0 \quad \text{and} \quad v = \log u \quad \text{in } \mathbf{R}^n \backslash \{0\}$$

is often used. Then, one should arrive at the following identities

$$\nabla u = e^v \nabla v = u \nabla v \tag{2.1}$$

and

$$\Delta u = u(|\nabla v|^2 + \Delta v) \tag{2.2}$$

and

$$\Delta^2 u = e^v \big(2|\nabla v|^2 \Delta v + 2\langle \nabla \Delta v, \nabla v\rangle + (\Delta v)^2 + \Delta^2 v + |\nabla v|^4 + 4\nabla^2 v(\nabla v, \nabla v) + \Delta |\nabla v|^2\big), \tag{2.3}$$

where $\nabla^2 v = (v_{ij})$ denotes the Hessian matrix of $v$. Making use of the above three identities we quickly obtain the following lemma whose proof can be found in [NN26b].

**Lemma 2.1.** *We have*

$$\begin{aligned} -\Delta^2 v = {} & 2|\nabla v|^2 \Delta v + 4\langle \nabla \Delta v, \nabla v\rangle + (\Delta v)^2 + |\nabla v|^4 \\ & + 4\nabla^2 v(\nabla v, \nabla v) + 2\|\nabla^2 v\|^2 - |x|^\sigma e^{2kv}, \end{aligned} \tag{2.4}$$

*where* $\|\nabla^2 v\|$ *denotes the Hilbert–Schmidt norm on matrices defined to be*

$$\|\nabla^2 v\| = \Big(\sum_{i,j} |v_{ij}|^2\Big)^{1/2}.$$

Thanks to (2.1), we know that

$$\alpha \frac{|\nabla u|^2}{u} + \beta |x|^{\frac{\sigma}{2}} u^{\frac{p+1}{2}} = u\big(\alpha |\nabla v|^2 + \beta |x|^{\frac{\sigma}{2}} e^{kv}\big) \quad \text{in } \mathbf{R}^n \backslash \{0\},$$

where $\alpha$ and $\beta$ are as follows

$$\alpha := \frac{2}{n-4} \quad \text{and} \quad \beta := \sqrt{\frac{2}{p+1-\frac{8}{n(n-4)}}}. \tag{2.5}$$

This and (2.2) allow to transform the desired estimate for $u$ into

$$0 \geqslant u\big(\Delta v + (\alpha+1)|\nabla v|^2 + \beta |x|^{\frac{\sigma}{2}} e^{kv}\big) \quad \text{in } \mathbf{R}^n \backslash \{0\}.$$

This suggests us to define the following auxiliary function $w : \mathbf{R}^n \backslash \{0\} \to \mathbf{R}$ given by

$$w = \Delta v + (\alpha+1)|\nabla v|^2 + \beta |x|^{\frac{\sigma}{2}} e^{kv}. \tag{2.6}$$

Our goal is to show that $w \leqslant 0$ everywhere in $\mathbf{R}^n \backslash \{0\}$. To this end, we first mention the following weaker version.

**Lemma 2.2.** *We have*

$$w \leqslant \beta |x|^{\frac{\sigma}{2}} e^{kv} \quad \textit{in } \mathbf{R}^n \backslash \{0\}. \tag{2.7}$$

Since the estimate (2.7) is equivalent to

$$u\Delta u + \alpha |\nabla u|^2 \leqslant 0 \quad \text{in } \mathbf{R}^n \backslash \{0\},$$

which is nothing but (1.4) in the punctured space $\mathbf{R}^n \backslash \{0\}$, which is now well-known, see [**NN26c**]. Hence, no proof is actually required.

Next, we compute and estimate $\Delta w$. This step is very important in Bernstein's method. That said, we have the following claim.

**Lemma 2.3.** *We have*

$$\begin{aligned}\Delta w \geqslant{}& 2(\alpha-1)\langle \nabla v, \nabla w\rangle - I_0 w^2 + I_1 \alpha w |\nabla v|^2 + I_4 \beta |x|^{\frac{\sigma}{2}} w e^{kv} \\ &+ I_5 \beta |x|^{\frac{\sigma}{2}} e^{kv} |\nabla v|^2 + \frac{\beta\sigma}{2}\Big(\big(n+\frac{\sigma}{2}-2\big) + 2(k-\alpha+1)\langle x, \nabla v\rangle\Big)|x|^{\sigma/2-2} e^{kv}\end{aligned} \tag{2.8}$$

*in* $\mathbf{R}^n \backslash \{0\}$, *where* $I_0$, $I_1$, $I_4$, *and* $I_5$ *are constant given in as follows*

$$\begin{cases} I_0 = 1 - \dfrac{2\alpha}{n}, \\ I_1 = \dfrac{2(n-6)}{n-4}, \\ I_4 = \dfrac{1}{\beta^2} + I_0, \\ I_5 = \dfrac{1}{4}\Big(p - \dfrac{n+4}{n-4}\Big)\Big(p + \dfrac{n-8}{n-4}\Big) \geqslant 0. \end{cases} \tag{2.9}$$

Proof of (2.8) is rather involved, but can be found in [NN26b]. From (2.9), it is worth noting that $I_0 > 0$ and $I_5 \geqslant 0$ due to the condition $p \geqslant (n+4+\sigma)/(n-4)$. See also Remark 3.2.

In the estimate (2.8) above, there is a term involving $\langle x, \nabla v\rangle$ which requires extra care. It is clear that this term appears mainly because of the weight function $|x|^\sigma$. This is also the reason why we have to keep the term involving $I_5$ in the estimate of $\Delta w$ although we know that $I_5 \geqslant 0$. To go further, let us introduce a threshold $\sigma^*$ for $\sigma$ defined by

$$\sigma^* = \frac{2(n-2)I_5}{(k-\alpha+1)^2 - I_5}.$$

This will play an important role, as we soon see. But first we show that $\sigma^*$ is well-defined.

**Lemma 2.4.** *We have*

$$(k-\alpha+1)^2 - I_5 > 0.$$

*In particular, there holds $\sigma^* \geqslant 0$,*

We now describe the role of the threshold $\sigma^*$. In fact, we can further estimate $\Delta w$ from the below provided $0 < \sigma \leqslant \sigma^*$.

**Lemma 2.5.** *Suppose $0 \leqslant \sigma \leqslant \sigma^*$. Then, we have*

$$\Delta w \geqslant 2(\alpha-1)\langle \nabla v, \nabla w\rangle - I_0 w^2 + I_1 \alpha w |\nabla v|^2 + I_4 \beta |x|^{\frac{\sigma}{2}} w e^{kv} \tag{2.10}$$

*in $\mathbf{R}^n \backslash \{0\}$.*

Having Lemma 2.5 in hand, to obtain the crucial estimate (2.10) for $\Delta w$ from the below, it suffices to examine the inequality $\sigma \leqslant \sigma^*$, Surprisingly, we show that this is always the case.

**Proposition 2.6.** *Under the conditions $n \geqslant 5$, $\sigma \geqslant 0$, and $p \geqslant (n+4+\sigma)/(n-4)$, there always holds $\sigma^* \geqslant \sigma$. Consequently, the key estimate (2.10) holds in $\mathbf{R}^n \backslash \{0\}$ for any $\sigma \geqslant 0$.*

*Proof.* This follows from direct verification. Indeed, by the definition of $\sigma^*$ we have

$$\sigma^* = \frac{2(n-2)I_5}{(k-\alpha+1)^2 - I_5} = \frac{\frac{1}{2}(n-2)(p - \frac{n+4}{n-4})(p + \frac{n-8}{n-4})}{\frac{1}{2}(p+\frac{n-8}{n-4})\frac{n-2}{n-4}} = \frac{p - \frac{n+4}{n-4}}{\frac{1}{n-4}} \geqslant \sigma,$$

where the latter estimate is exactly the the inequality $p \geqslant (n+4+\sigma)/(n-4)$, which was used for the second time. From this we obtain the key estimate (2.10) for all $\sigma \geqslant 0$. □

Our final piece of preparation is a useful pointwise estimate. Before proving such an estimate, we recall that any solution $u$ to $(1.1)_\sigma$ in $\mathbf{R}^n \backslash \{0\}$ with $n \geqslant 5$, $\sigma > -4$, and $p > (n+\sigma)/(n-4)$ enjoys the following integral representation

$$u(x) = C(4) \int_{\mathbf{R}^n} \frac{|y|^\sigma u^p(y)}{|x-y|^{n-4}} dy \quad \text{in } \mathbf{R}^n \backslash \{0\} \tag{2.11}$$

and

$$\nabla u(x) = -(n-4)C(4) \int_{\mathbf{R}^n} \frac{x-y}{|x-y|^{n-2}} |y|^\sigma u^p(y) dy \quad \text{in } \mathbf{R}^n \backslash \{0\} \tag{2.12}$$

and

$$-\Delta u(x) = C(2)\int_{\mathbf{R}^n} \frac{|y|^\sigma u^p(y)}{|x-y|^{n-2}}dy \quad \text{in } \mathbf{R}^n\backslash\{0\} \tag{2.13}$$

with

$$C(4) = \frac{\Gamma(\frac{n-4}{2})}{16\pi^{n/2}}, \quad C(2) = \frac{\Gamma(\frac{n-2}{2})}{4\pi^{n/2}} = 2(n-4)C(4). \tag{2.14}$$

See [**CAM08**] and [**NY22**]. Then we can prove the following.

**Lemma 2.7.** *Assume $n \geqslant 5$, $p > (n+\sigma)/(n-4)$ and $\sigma > -4$. Then, any radial solution $u$ is monotone deceasing in $(0,+\infty)$ and there exist a positive constant $C_0$ such that*

$$u(r) \leqslant C_0 r^{-\theta}, \quad 0 \leqslant |\nabla u|(r) \leqslant C_0 r^{-1} u(r), \quad |\Delta u|(r) \leqslant C_0 r^{-\theta-2}$$

*in $(0,+\infty)$.*

This is more or less well-known, and a proof can be found in [**NN26c**]. Using the integral equation (2.11), we know that any solution, not necessarily radial, has the following lower bound, for e.g. see [**NY22**].

**Lemma 2.8.** *Assume $n \geqslant 5$, $p > (n+\sigma)/(n-4)$ and $\sigma > -4$. There is some constant $C > 0$ such that*

$$u(x) \geqslant \frac{C}{|x|^{n-4}} \quad \textit{for } |x| \geqslant 2.$$

As a consequence of above estimates, we have the following crucial estimate for $w$ from the above.

**Lemma 2.9.** *Assume $n \geqslant 5$, $p > (n+\sigma)/(n-4)$ and $\sigma > -4$. Then, there exist a positive (uniform) constant $C_0$ such that*

$$w(r) \leqslant C_0 r^{-2} \quad \textit{in } (0,+\infty).$$

*Proof.* From (2.13) and the definition of $w$ in (2.6) we have

$$w = \frac{\Delta u}{u} + \alpha\Big(\frac{|\nabla u|}{u}\Big)^2 + \beta r^{\sigma/2} u^{(p-1)/2} < \alpha\Big(\frac{|\nabla u|}{u}\Big)^2 + \beta r^{\sigma/2} u^{(p-1)/2}$$

for some $C_0 > 0$. From this and Lemma 2.7 we easily obtain the desired estimate. □

It would be interesting to know if $w$ can be estimated from the below. For e.g. should we have

$$|w(r)| \leqslant C_0 r^{-2} \quad \text{in } (0,+\infty)$$

instead of the one-sided estimate in Lemma 2.9 above?

Our last two lemmas are gradient estimates. We start with a gradient estimate for $u$, which is simply a consequence of (1.4).

**Lemma 2.10.** *Assume $n \geqslant 5$, $p > (n+\sigma)/(n-4)$ and $\sigma > -4$. There holds*

$$n-4 \geqslant -\frac{ru'(r)}{u(r)} > 0 \quad \textit{in } (0,+\infty). \tag{2.15}$$

*Proof.* This can be proved by using an idea from [**CMM93**], as done in [**NN26b**]. Here we offer a new proof by exploiting the role of the inequality (1.4). We first denote

$$v = \log u.$$

Then, in terms of $v$, the differential inequality (1.4) becomes

$$-\Delta v \geqslant \frac{n-2}{n-4}|\nabla v|^2 \quad \text{in } \mathbf{R}^n\backslash\{0\},$$

which then yields

$$0 \geqslant r^2 v'' + (n-1)rv' + \frac{n-2}{n-4}r^2(v')^2 \quad \text{in } (0,+\infty).$$

(Here note that $v$ is also radial.) For simplicity we denote

$$\phi = rv' = r\frac{u'}{u}.$$

Then the preceding inequality becomes

$$0 \geqslant r\phi' + (n-2)\phi + \frac{n-2}{n-4}\phi^2 = r\phi' + \frac{n-2}{n-4}\phi(\phi+n-4)$$

in $(0,+\infty)$. Our aim is to show that

$$\phi + n - 4 \geqslant 0 \quad \text{in } (0,+\infty).$$

By way of contradiction suppose that there is some $r_0 > 0$ such that

$$\phi(r_0) + n - 4 < 0.$$

Keep in mind that $\phi < 0$ because $u$ is monotone decreasing. Then

$$0 > -\frac{n-2}{n-4}\phi(\phi+n-4) \geqslant r\phi' \quad \text{at } r_0.$$

Hence, $\phi$ is decreasing right after $r_0$, namely $\phi + n - 4 < 0$ right after $r_0$. The smoothness of $\phi$ tells us that we actually have

$$\phi + n - 4 < 0 \quad \text{for all } r \geqslant r_0.$$

We now denote

$$C = \frac{n-2}{n-4}\phi(r_0)(\phi(r_0)+n-4) > 0.$$

Then,

$$\phi'(r) \leqslant -\frac{C}{r} \quad \text{for all } r \geqslant r_0.$$

Integrating both sides give

$$\phi(r) \leqslant \phi(r_0) - C\log\frac{r}{r_0} < -(n-4) - C\log\frac{r}{r_0} \quad \text{for all } r \geqslant r_0.$$

Thus,

$$r(\log u)' < -(n-4+C) \quad \text{for all } r \geqslant 3r_0,$$

namely

$$\log u < \log(3r_0u(3r_0)) - (n-4+C)\log r \quad \text{for all } r \geqslant 3r_0.$$

Thus, we have just shown that

$$u(r) \leqslant \frac{3r_0u(3r_0)}{r^{n-4+C}} \quad \text{for all } r \geqslant 3r_0,$$

which violates the lower bound of $u$ in Lemma 2.8, thanks to $C > 0$. This completes the present proof. □

Our last lemma in this section is a gradient estimate for $\Delta u$.

**Lemma 2.11.** *Assume $n \geqslant 5$, $p > (n+\sigma)/(n-4)$ and $\sigma > -4$. There holds*

$$n-2 \geqslant -\frac{r(\Delta u)'(r)}{\Delta u(r)} > 0 \quad \text{in } (0,+\infty). \tag{2.16}$$

*Proof.* Observe

$$(r(-\Delta u)'(r) + (n-2)(-\Delta u)(r))' = r\Big((-\Delta u)''(r) + \frac{n-1}{r}(-\Delta u)'(r)\Big) = -r\Delta^2 u < 0$$

in $(0,+\infty)$. Hence

$$r \mapsto r(-\Delta u)'(r) + (n-2)(-\Delta u)(r)$$

is strictly decreasing. If there were some $r_0 > 0$ such that

$$r_0(-\Delta u)'(r_0) + (n-2)(-\Delta u)(r_0) < 0,$$

then we would have

$$r(-\Delta u)'(r) + (n-2)(-\Delta u)(r) < r_0(-\Delta u)'(r_0) + (n-2)(-\Delta u)(r_0) =: C < 0$$

in $[r_0,+\infty)$, namely

$$(-\Delta u)'(r) < (-\Delta u)'(r) + \frac{n-2}{r}(-\Delta u)(r) < \frac{C}{r} < 0.$$

Integrating over $[r_0, r]$ gives

$$-\infty < -(-\Delta u)(r_0) < (-\Delta u)(r) - (-\Delta u)(r_0) < C \log \frac{r}{r_0}$$

in $[r_0,+\infty)$. Sending $r \nearrow +\infty$ we obtain contradiction, hence

$$r(-\Delta u)' + (n-2)(-\Delta u) \geqslant 0 \quad \text{in } (0,+\infty).$$

With help of [**NNT26**, Remark A.2(3)], the positivity of $-r(\Delta u)'(r)/\Delta u(r)$ in $(0,+\infty)$ is clear, thanks to the $\Delta u < 0$ everywhere in $\mathbf{R}^n \backslash \{0\}$, see (2.13). Putting the two estimates above gives (2.16). □

It would be interesting if Lemma 2.11 can be proved by using the method used in the proof of Lemma 2.10. Perhaps, looking for an inequality for $\Delta u$ in the same fashion of (1.4) for is needed.

## 3. Proof of Theorem 1.1 and remarks

Let $u$ be a non-negative, radial, $C^4$-solution to $(1.1)_\sigma$ in $\mathbf{R}^n \backslash \{0\}$ with $n \geqslant 5$, $\sigma \geqslant 0$, and $p \geqslant (n+4+\sigma)/(n-4)$. Then, it is well-known that $u > 0$, for e.g. see the integral representation (2.11).

3.1. **Proof of Theorem 1.1.** This section is devoted to the proof of Theorem 1.1. Since the sharpness of the coefficient $2/(n-4)$ follows from Theorem 1.3 and the optimality of the threshold $(n+4+\sigma)/(n-4)$ is mentioned in Remark 3.2, it suffices to prove (1.6) for radial solutions $u$ to $(1.1)_\sigma$ in $\mathbf{R}^n \backslash \{0\}$. In addition, since the standard singular solution $u_{\rm sing}$ always satisfies the inequality, even with better coefficients, we assume that $u$ is not the standard solution.

Depending on the size of $p$ we have three possible cases as follows.

**Case 1.** Suppose $(n+4+\sigma)/(n-4) \leqslant p < (n+4+2\sigma)/(n-4)$. Let $u$ be any radial

solution to $(1.1)_\sigma$ in $\mathbf{R}^n\backslash\{0\}$. We denote by $\widehat u$ the Kelvin transform of $u$, namely

$$\widehat u(x) = \frac{1}{|x|^{n-4}} u\Big(\frac{x}{|x|^2}\Big) \quad \text{in } \mathbf{R}^n\backslash\{0\}.$$

Then, we know that $\widehat u$ solves $(1.1)_{\widehat\sigma}$ in $\mathbf{R}^n\backslash\{0\}$, namely

$$\Delta^2 \widehat u = |x|^{\widehat\sigma} \widehat u^p \quad \text{in } \mathbf{R}^n\backslash\{0\} \tag{3.1}$$

with the following new conditions

$$\widehat\sigma = (n-4)p - (n+4+\sigma) \geqslant 0$$

and

$$p - \frac{n+4+2\widehat\sigma}{n-4} = \frac{n+4+2\sigma}{n-4} - p > 0$$

due to

$$\frac{n+4+\sigma}{n-4} \leqslant p < \frac{n+4+2\sigma}{n-4}.$$

Also, we have the radial symmetry of $\widehat u$ as well. Hence, the transformed equation (3.1) solved by $\widehat u$ is supercritical. Kêp in mind that $u$ is not the standard solution. Consequently, the radial symmetry of $\widehat u$ allows us to apply [**NN26c**, Theorem 1.2] to the transformed equation (3.1) to conclude that $\widehat u$ is regular at the origin, namely $\widehat u(0) > 0$ and $\widehat u \in C(\mathbf{R}^n) \cap C^4(\mathbf{R}^n\backslash\{0\})$ in the extension sense. Thanks to $\widehat\sigma \geqslant 0$, a general regularity result in [**GW24**] tells us that we even have $\widehat u \in C^4(\mathbf{R}^n)$ in the extension sense. Therefore, we can deduce from [**NN26b**] the validity of

$$-\Delta \widehat u(x) \geqslant \frac{2}{n-4}\frac{|\nabla \widehat u|^2}{\widehat u}(x) + \sqrt{\frac{2}{p+1-\frac{8}{n(n-4)}}} |x|^{\frac{\widehat\sigma}{2}} \widehat u^{\frac{p+1}{2}}$$

in $\mathbf{R}^n\backslash\{0\}$. Notice that

$$\Delta \widehat u(x) + \frac{2}{n-4}\frac{|\nabla \widehat u|^2}{\widehat u}(x) = \frac{1}{|x|^n}\Big[\Delta u\Big(\frac{x}{|x|^2}\Big) + \frac{2}{n-4}\frac{|\nabla u|^2}{u}\Big(\frac{x}{|x|^2}\Big)\Big] \tag{3.2}$$

and

$$|x|^{\frac{\widehat\sigma}{2}} \widehat u^{\frac{p+1}{2}}(x) = \frac{1}{|x|^n}\Big|\frac{x}{|x|^2}\Big|^{\frac{\sigma}{2}} u^{\frac{p+1}{2}}\Big(\frac{x}{|x|^2}\Big) \tag{3.3}$$

in $\mathbf{R}^n\backslash\{0\}$; see also [**NN26b**]. Therefore, we are in position to make use of (3.2) and (3.3) to go back from $\widehat u$ to $u$, with the scaling factor $|x|^{-n}$, to obtain the desired inequality (1.6).

**Case 2**. Suppose $p > (n+4+2\sigma)/(n-4)$. In this case, we don't need to make use of the Kelvin transform because the argument in Case 1 above works well with the transformed function $\widehat u$ replaced by the original function $u$. Hence, details are left for the readers, and this completes the proof in this case.

**Case 3**. Suppose $p = (n+4+2\sigma)/(n-4)$. In this case, using the classification of radial solutions it suffices to consider the case $|x|^{(n-4)/2}u(x)$ is bounded. As a consequence, it follows from Lemma 2.9 that there is some $C > 0$ such that

$$w(x) \leqslant \frac{C}{|x|^2} \quad \text{in } \mathbf{R}^n\backslash\{0\}.$$

(Here we only need the upper bound for $w$.) Let

$$\widehat w(x) = |x|^2 w(x) \quad \text{in } \mathbf{R}^n\backslash\{0\}, \tag{3.4}$$

or equivalently,

$$w = |x|^{-2}\widehat w \quad \text{in } \mathbf{R}^n\backslash\{0\}.$$

Then, we quickly have

$$\nabla w = |x|^{-2}\nabla \widehat{w} - 2|x|^{-4} x\, \widehat{w}$$

and

$$\nabla v \cdot \nabla w = |x|^{-2}\nabla v \cdot \nabla \widehat{w} - 2|x|^{-4}(x\cdot \nabla v)\widehat{w}$$

in $\mathbf{R}^n \backslash \{0\}$. Moreover,

$$\Delta w = |x|^{-2}\Big(\Delta \widehat{w} - \frac{4}{|x|^2} x\cdot \nabla \widehat{w} - \frac{2(n-4)}{|x|^2}\widehat{w}\Big) \quad \text{in } \mathbf{R}^n \backslash \{0\}.$$

Substituting these identities into the key estimate (2.10), namely

$$\Delta w \geqslant 2(\alpha-1)\langle \nabla v, \nabla w\rangle - I_0 w^2 + I_1 \alpha w|\nabla v|^2 + I_4 \beta |x|^{\frac{\sigma}{2}} w e^{kv} \quad \text{in } \mathbf{R}^n \backslash \{0\},$$

and multiplying both sides of the resulting inequality by $|x|^2$, we arrive at

$$\begin{aligned}\Delta \widehat{w} - \frac{4}{|x|^2} x\cdot \nabla \widehat{w} - \frac{2(n-4)}{|x|^2}\widehat{w} \geqslant{}& 2(\alpha-1)\nabla v\cdot \nabla \widehat{w} - \frac{4(\alpha-1)}{|x|^2}(x\cdot \nabla v)\widehat{w}\\ &- I_0 |x|^{-2}\widehat{w}^2 + I_1 \alpha |\nabla v|^2 \widehat{w} + I_4 \beta |x|^{\frac{\sigma}{2}} e^{kv}\widehat{w}\end{aligned}$$

in $\mathbf{R}^n \backslash \{0\}$. Equivalently, we should have

$$\begin{aligned}\Delta \widehat{w} - \frac{4}{|x|^2} x\cdot \nabla \widehat{w} - 2(\alpha-1)\nabla v\cdot \nabla \widehat{w} \geqslant{}& \Big(\frac{2(n-4)}{|x|^2} - \frac{4(\alpha-1)}{|x|^2} x\cdot \nabla v + I_1 \alpha |\nabla v|^2\Big)\widehat{w}\\ &- I_0 |x|^{-2}\widehat{w}^2 + I_4 \beta |x|^{\frac{\sigma}{2}} e^{kv}\widehat{w}.\end{aligned}$$

in $\mathbf{R}^n \backslash \{0\}$. Since $u$ is radial, so are $v$, $w$, and $\widehat{w}$. Difference from our previous work [NN26b], our function $\widehat{w}$ does not define at the origin, hence the standard cut-off function does not work since we have an obstacle at the origin. To avoid this, we use the so-called Emden–Fowler variable

$$s = \log |x| \quad \text{in } \mathbf{R}.$$

The idea is to work in one dimensional space $\mathbf{R}$ instead of $\mathbf{R}^n \backslash \{0\}$. Making use of the Emden–Fowler variable we write

$$\widehat{w}(x) = \widehat{W}(s), \quad v(x) = V(s). \tag{3.5}$$

Since

$$\frac{d}{dr} = \frac{1}{r}\frac{d}{ds},$$

we have

$$\widehat{w}_r = \frac{1}{r}\widehat{W}_s, \quad \widehat{w}_{rr} = \frac{1}{r^2}(\widehat{W}_{ss} - \widehat{W}_s),$$

and therefore

$$\Delta \widehat{w} = \widehat{w}_{rr} + \frac{n-1}{r}\widehat{w}_r = \frac{1}{r^2}\big(\widehat{W}_{ss} + (n-2)\widehat{W}_s\big).$$

Moreover,

$$x\cdot \nabla \widehat{w} = r\widehat{w}_r = \widehat{W}_s,$$

and

$$\nabla v\cdot \nabla \widehat{w} = v_r \widehat{w}_r = \frac{1}{r^2} V_s \widehat{W}_s.$$

Since

$$|\nabla v|^2 = \frac{1}{r^2} V_s^2,$$

and

$$|x|^{\frac{\sigma}{2}} e^{kv} = e^{\frac{\sigma}{2}s + kV},$$

the differential inequality becomes

$$\frac{1}{r^2}\big(\widehat{W}_{ss}+(n-2)\widehat{W}_s\big)-\frac{4}{r^2}\widehat{W}_s-\frac{2(\alpha-1)}{r^2}V_s\widehat{W}_s$$
$$\geqslant\Big(\frac{2(n-4)}{r^2}-\frac{4(\alpha-1)}{r^2}V_s+\frac{I_1\alpha}{r^2}V_s^2\Big)\widehat{W}-\frac{I_0}{r^2}\widehat{W}^2+I_4\beta e^{\frac{\sigma}{2}s+kV}\widehat{W}$$

in $\mathbf{R}$. Multiplying both sides by $r^2=e^{2s}$ yields

$$\begin{aligned}\widehat{W}_{ss}+(n-6)\,\widehat{W}_s-2(\alpha-1)V_s\widehat{W}_s\geqslant&\big(2(n-4)-4(\alpha-1)V_s+I_1\alpha V_s^2\big)\widehat{W}\\&-I_0\widehat{W}^2+I_4\beta e^{(\frac{\sigma}{2}+2)s+kV}\widehat{W}\end{aligned}$$

in $\mathbf{R}$. Finally, introducing the Emden–Fowler transform

$$Z(s)=e^{\frac{n-4}{2}s}u(e^s)\quad\text{in }\mathbf{R},$$

we have

$$V(s)=\log Z(s)-\frac{n-4}{2}s$$

and

$$V_s(s)=\frac{Z_s(s)}{Z(s)}-\frac{n-4}{2}$$

in $\mathbf{R}$. Keep in mind that $k=(p-1)/2$. Hence, there holds

$$e^{(\frac{\sigma}{2}+2)s+kV}=Z^{\frac{p-1}{2}},$$

and therefore

$$\begin{aligned}\widehat{W}_{ss}+(n-6)\,\widehat{W}_s-2(\alpha-1)V_s\widehat{W}_s\geqslant&\big(2(n-4)-4(\alpha-1)V_s+I_1\alpha V_s^2\big)\widehat{W}\\&-I_0\widehat{W}^2+I_4\beta Z^{\frac{p-1}{2}}\widehat{W}.\end{aligned}\tag{3.6}$$

in $\mathbf{R}$. So far we have seen that to avoid the singularity at the origin, we use the Emden–Fowler variable, in other words, we transfer the problem in the punctured space to the problem in one dimensional space. Consequently, we are forced to use one dimensional cut-off functions. By way of contradiction we suppose that

$$\widehat{W}(s_*)=\delta>0\tag{3.7}$$

for some $s_*\in\mathbf{R}$. To this purpose, let $\eta\in C_c^\infty(\mathbb{R})$ satisfy

$$0\leqslant\eta\leqslant 1,\quad \eta\equiv 1\text{ on }[-1,1],\quad \eta\equiv 0\text{ outside }[-2,2].$$

For $R>\max\{|s_*|,1\}$, we define

$$\eta_R(s)=\eta\big(\frac{s}{R}\big).$$

Then, there is some $C>0$ such that

$$|\eta_R'|\leqslant\frac{C}{R},\quad |\eta_R''|\leqslant\frac{C}{R^2}\quad\text{in }\mathbf{R}.\tag{3.8}$$

We now define

$$\Phi_R(s)=\eta_R(s)\widehat{W}(s)\quad\text{in }\mathbf{R}.$$

Notice that $\eta_R(s_*)=1$, which then yields

$$\Phi_R(s_*)=\widehat{W}(s_*)=\delta.$$

Now observe that $\Phi_R\equiv 0$ outside $[-2R,2R]$. Therefore, there is some maximum point of $\Phi_R$ in $\mathbf{R}$, denoted by $s_R\in\mathbb{R}$, namely

$$\Phi_R(s_R)=\max_{\mathbb{R}}\Phi_R.$$

Then by the contradiction argument, we have that

$$\Phi_R(s_R)\geqslant\Phi_R(s_*)=\delta>0.$$

In particular, there holds

$$\widehat{W}(s_R) > 0, \quad \eta_R(s_R) > 0.$$

In view of (3.5), (3.4), and Lemma 2.9 we know that

$$\widehat{W} \leqslant C_0 \quad \text{in } \mathbf{R},$$

where $C_0$ is the constant given in Lemma 2.9. Hence,

$$\eta_R(s_R) \geqslant \frac{\delta}{C_0}$$

which implies

$$\frac{s_R}{R} \in \Big\{ t \in \mathbf{R} : \eta(r) \geqslant \frac{\delta}{C_0} \Big\} =: K \subset [-R, R].$$

Here $K$ is a compact set. Note that the two numbers $\delta$ and $C_0$ as well as the set $K$ do not depend on $R$. Thus,

$$\sup_K \Big| \frac{\eta'}{\eta} \Big| < +\infty, \quad \sup_K \Big| \frac{\eta''}{\eta} \Big| < +\infty.$$

Note that $s_R/R \in K$. Hence, for some possibly new constant $C > 0$, see (3.8), we have

$$\Big| \frac{\eta_R'}{\eta_R} \Big| \leqslant \frac{C}{R}, \quad \Big| \frac{\eta_R''}{\eta_R} \Big| \leqslant \frac{C}{R^2} \quad \text{at } s_R. \tag{3.9}$$

(Compared with (3.8), the estimate (3.9) is possibly not true in the compact support of $\eta_R$; however, we only need it at $s_R$, which is true.) At the maximum point $s_R$, there hold

$$\Phi_R'(s_R) = 0, \quad \Phi_R''(s_R) \leqslant 0.$$

Solving the identity $\Phi_R' = \eta_R' \widehat{W} + \eta_R \widehat{W}_s$ gives

$$\widehat{W}_s = -\frac{\eta_R'}{\eta_R} \widehat{W} \quad \text{at } s_R.$$

Moreover, as

$$\Phi_R'' = \eta_R'' \widehat{W} + 2\eta_R' \widehat{W}_s + \eta_R \widehat{W}_{ss}$$

and $\eta_R(s_R) > 0$, we obtain

$$0 \geqslant \frac{\eta_R''}{\eta_R} \widehat{W} + \frac{2\eta_R'}{\eta_R} \widehat{W}_s + \widehat{W}_{ss} \quad \text{at } s_R.$$

Combining with the expression of $\widehat{W}_s$ and the estimate (3.9) yields

$$\begin{aligned}
0 &\geqslant \widehat{W}_{ss} + \Big( \frac{\eta_R''}{\eta_R} - 2\frac{(\eta_R')^2}{\eta_R^2} \Big) \widehat{W}, \\
&\geqslant \widehat{W}_{ss} - \frac{C}{R^2} \widehat{W} \\
&\geqslant -\frac{C}{R^2} \widehat{W} - (n-6)\widehat{W}_s + 2(\alpha-1)V_s \widehat{W}_s + \big(2(n-4) - 4(\alpha-1)V_s + I_1 \alpha V_s^2\big) \widehat{W} \\
&\quad - I_0 \widehat{W}^2 + I_4 \beta Z^{\frac{p-1}{2}} \widehat{W} \\
&= \Big( -\frac{C}{R^2} + (n-6)\frac{\eta_R'}{\eta_R} - 2(\alpha-1)V_s \frac{\eta_R'}{\eta_R} + 2(n-4) - 4(\alpha-1)V_s + I_1 \alpha V_s^2 \Big) \widehat{W} \\
&\quad - I_0 \widehat{W}^2 + I_4 \beta Z^{\frac{p-1}{2}} \widehat{W}
\end{aligned}$$

at $s_R$. (Here $C > 0$ is a generic constant, which could change from line to line.) We now estimate the last two terms on the right side of the preceding long formula. Clearly,

$$-I_0\widehat{W}^2 + I_4\beta Z^{\frac{p-1}{2}}\widehat{W} = I_4\beta|x|^{\frac{\sigma}{2}+2}e^{kv}\widehat{w} - I_0\widehat{w}^2 \geqslant (I_4 - I_0)\widehat{w}^2 = \beta^{-2}\widehat{w}^2, \tag{3.10}$$

see (2.9). So

$$-I_0\widehat{W}^2 + I_4\beta Z^{\frac{p-1}{2}}\widehat{W} \geqslant \beta^{-2}\widehat{W}^2 \geqslant 0 \quad \text{in } \mathbf{R}.$$

Now depending on the size of $n$ we have two subcases.

*Subcase 3.1.* Suppose $n \geqslant 6$. Then we quickly have

$$\begin{aligned} 2(n-4) - 4(\alpha-1)V_s + I_1\alpha V_s^2 &= 2(n-4) + \frac{4(n-6)}{(n-4)}V_s + \frac{4(n-6)}{(n-4)^2}V_s^2 \\ &= \frac{4(n-6)}{(n-4)^2}\big(V_s + \frac{n-4}{2}\big)^2 + n - 2 \geqslant n-2 \end{aligned}$$

everywhere in $\mathbf{R}$. Then, we arrive at

$$0 \geqslant \Big(-\frac{C}{R^2} + (n-6)\frac{\eta_R'}{\eta_R} - 2(\alpha-1)V_s\frac{\eta_R'}{\eta_R} + n - 2\Big)\widehat{W}$$

at $s_R$. Keep in mind that by definition $\widehat{W}(s_R) > 0$ and is bounded from above. Now, making use of (3.8)–(3.9) and the boundedness of $V_s$ as seen from the identity

$$|V_s| = \frac{r|\nabla u|}{u}$$

we deduce that

$$-\frac{C}{R^2} + (n-6)\frac{\eta_R'(s_R)}{\eta_R(s_R)} - 2(\alpha-1)V_s(s_R)\frac{\eta_R'(s_R)}{\eta_R(s_R)} \to 0 \tag{3.11}$$

as $R \to \infty$. Consequently, for large $R \gg 1$ we should have

$$0 \geqslant (n-3)\widehat{W}(s_R) > 0,$$

which is impossible. In view of (3.7), the contradiction shows that $\widehat{W}(s) \leqslant 0$ for all $s \in \mathbf{R}$, so is $\widehat{w}(x) \leqslant 0$ for all $x \in \mathbf{R}^n \backslash \{0\}$. Thus, $w(x) \leqslant 0$ in $\mathbf{R}^n \backslash \{0\}$ as claimed.

*Subcase 3.2.* Suppose $n = 5$. In this case we have

$$0 \geqslant \Big(-\frac{C}{R^2} + (n-6)\frac{\eta_R'}{\eta_R} - 2(\alpha-1)V_s\frac{\eta_R'}{\eta_R} + \frac{4(n-6)}{(n-4)^2}\big(V_s + \frac{n-4}{2}\big)^2 + n - 2\Big)\widehat{W},$$

which becomes

$$0 \geqslant \Big(-\frac{C}{R^2} - \frac{\eta_R'}{\eta_R} - 2(\alpha-1)V_s\frac{\eta_R'}{\eta_R} - 4\big(V_s + \frac{1}{2}\big)^2 + 3\Big)\widehat{W}.$$

We are now in position to make use of Lemma 2.10 to get

$$0 \geqslant V_s = \frac{ru'}{u} \geqslant -n + 4.$$

Hence, under $n = 5$, we have $|V_s + 1/2| \leqslant 1/2$, which further implies

$$0 \geqslant \Big(-\frac{C}{R^2} - \frac{\eta_R'}{\eta_R} - 2(\alpha-1)V_s\frac{\eta_R'}{\eta_R} + 2\Big)\widehat{W}$$

at $s_R$. Making use of (3.11) we have

$$0 \geqslant \widehat{W}(s_R) > 0$$

provided $R \gg 1$, which is impossible. The contradiction shows that $\widehat{W}(s) \leqslant 0$ for all $s \in \mathbf{R}$.

3.2. **Remarks on the thresholds of $\sigma$ and $p$.** As mentioned in Introduction, the following two remarks concern the validity of (1.6) with either $\sigma < 0$ or $p < (n+4+\sigma)/(n-4)$. We should remind that under $n \geqslant 5$ and $\sigma > -4$, either $1 < p < (n+2+\sigma)/(n-2)$ or $p > (n+\sigma)/(n-4)$ is necessarily for the existence of solutions to $(1.1)_\sigma$ in $\mathbf{R}^n \backslash \{0\}$, see [Ngo24]. However, by testing the left hand side of (1.6) by the standard singular solution of the form $L|x|^{-\theta}$ for some $L > 0$ with $\theta = (4+\sigma)/(p-1)$ we notice that

$$-\Delta(L|x|^{-\theta}) < 0 \quad \text{in } \mathbf{R}^n \backslash \{0\}$$

provided $1 < p < (n+2+\sigma)/(n-2)$. Hence, for the validity of (1.6), it is necessary to assume $p > (n+\sigma)/(n-4)$.

Now we have a remark concerning $-2 < \sigma < 0$, see Table 1.

*Remark* 3.1. Suppose that $n \geqslant 5$, $p \geqslant (n+4+2\sigma)/(n-4)$, and $\sigma \in (-2,0)$. Then, thanks to an existence result in [HS21], see also [NY22] for a generalization, there is some radial solution $u$ to $(1.1)_\sigma$ in $\mathbf{R}^n \backslash \{0\}$ satisfying

$$u \in C(\mathbf{R}^n) \cap C^4(\mathbf{R}^n \backslash \{0\}).$$

Therefore, by a regularity result in [GW24], we deduce that $u \in C^2(\mathbf{R}^n)$. Hence, near the origin, this and $\sigma < 0$ make two sides of (1.8) irrelevant.

Next, we have a remark concerning the regime $(n+2+\sigma)/(n-4) < p < (n+4+\sigma)/(n-4)$ and $\sigma > -2$, see also Table 1.

*Remark* 3.2. Let $u$ be a solution to $(1.1)_\sigma$ in $\mathbf{R}^n \backslash \{0\}$ with $\sigma > -2$ and $(n+2+\sigma)/(n-4) < p < (n+4+\sigma)/(n-4)$, which will be specified later. If $(n+4+2\sigma)/(n-4) \leqslant p < (n+4+\sigma)/(n-4)$, implicitly giving $\sigma < 0$, then the conclusion follows from Remark 3.1 above. Now suppose $(n+2+\sigma)/(n-4) < p < (n+4+2\sigma)/(n-4)$ and any $\sigma > -2$. Still by the Kelvin transform $\widehat{u}$, we again arrive at the equation (3.1)

$$\Delta^2 \widehat{u} = |x|^{\widehat{\sigma}} \widehat{u}^p \quad \text{in } \mathbf{R}^n \backslash \{0\}$$

solved by $\widehat{u}$ with the following new conditions

$$-4 < \widehat{\sigma} < 0 \quad \text{and} \quad p > \frac{n+4+2\widehat{\sigma}}{n-4}.$$

From this we can choose $\widehat{u}$ as a radial solution to $(1.1)_\sigma$ in $\mathbf{R}^n \backslash \{0\}$ in such a way that $\widehat{u} \in C(\mathbf{R}^n) \cap C^4(\mathbf{R}^n \backslash \{0\})$, see [HS21] and [NY22]. (From this, the desired solution $u$ will be found by taking the Kelvin transform of $\widehat{u}$.) Hence, if we further require $\widehat{\sigma} \in (-2,0)$, equivalently

$$\frac{n+2+\sigma}{n-4} < p < \frac{n+4+\sigma}{n-4},$$

then there holds $\widehat{u} \in C^2(\mathbf{R}^n)$ by the regularity result in [GW24]. Now the pointwise inequality (1.8) for $u$ becomes

$$-\Delta \widehat{u} \geqslant \frac{2}{n-4} \frac{|\nabla \widehat{u}|^2}{\widehat{u}} + C|x|^{\frac{\widehat{\sigma}}{2}} \widehat{u}^{\frac{p+1}{2}} \quad \text{in } \mathbf{R}^n \backslash \{0\},$$

see (3.2)–(3.3). With the above regularity in hand, the preceding inequality immediately provides contradiction due to $\widehat{\sigma} < 0$.

## 4. Proof of Theorem 1.3

As mentioned in Introduction, the pointwise inequality can be proved by either using integral representation [**LX25**] or by using the method of invariant tensors [**MWW25**]. Here we argue by using integral representation, thanks to (2.11)–(2.13).

Indeed, let $u \in C^4(\mathbf{R}^n\backslash\{0\})$ be a solution to $(1.1)_\sigma$, and we recall from (2.11) the integral representation

$$\nabla u(x) = -(n-4)C(4)\int_{\mathbf{R}^n} \frac{x-y}{|x-y|^{n-2}}|y|^\sigma u^p(y)dy \quad \text{in } \mathbf{R}^n\backslash\{0\}.$$

In particular, $u>0$ and using (2.12) there holds

$$|\nabla u|(x) \leqslant (n-4)C(4)\int_{\mathbf{R}^n} \frac{|y|^\sigma u^p(y)}{|x-y|^{n-3}}dy \quad \text{in } \mathbf{R}^n\backslash\{0\}.$$

Hence, by Hölder's inequality and (2.13) we eventually have

$$\begin{aligned}\frac{2}{n-4}|\nabla u|^2 &\leqslant 2(n-4)C(4)^2\Big(\int_{\mathbf{R}^n} \frac{|y|^\sigma u^p(y)}{|x-y|^{n-3}}dy\Big)^2\\ &\leqslant C(2)C(4)\Big(\int_{\mathbf{R}^n} \frac{|y|^\sigma u^p(y)}{|x-y|^{n-4}}dy\Big)\Big(\int_{\mathbf{R}^n} \frac{|y|^\sigma u^p(y)}{|x-y|^{n-2}}dy\Big) = -u\Delta u\end{aligned}$$

in $\mathbf{R}^n\backslash\{0\}$. Hence, in the rest of the present proof, we show the sharpness of the constant $2/(n-4)$.

For clarity, we denote by $\alpha_{\rm GHWW}$ the sharp constant, which satisfies $\alpha_{\rm GHWW} \geqslant 2/(n-4)$. Then, it suffices to prove that $\alpha_{\rm GHWW} = 2/(n-4)$. To this purpose we first let $U$ be a radial solution to $(1.1)_\sigma$ in $\mathbf{R}^n\backslash\{0\}$, which is not the standard singular solution $L|x|^{-\theta}$. Such an existence can be found in [**NY22**, Corollary 7.6] for the subcritical case $(n+\sigma)/(n-4) < p < (n+4+2\sigma)/(n-4)$, in [**JL14**, Theorem 3] for the critical case $p=(n+4+2\sigma)/(n-4)$, and again in [**NY22**, Corollary 7.6] for the supercritical case $p>(n+4+2\sigma)/(n-4)$ after making use of the Kelvin transform. Clearly, the radial solution $U$ enjoys the following asymptotic behavior at infinity

$$U(|x|) \sim \frac{1}{|x|^{n-4}} \quad \text{for } |x| \gg 1. \tag{4.1}$$

Here by $f\sim g$ we mean that $|f(x)/g(x)|$ is bounded near infinity.

The estimate (4.1) for any $p>(n+\sigma)/(n-4)$ requires some care as we are about to explain. In the subcritical case $(n+\sigma)/(n-4) < p < (n+4+2\sigma)/(n-4)$, the estimate (4.1) is guaranteed from [**NY22**, Corollary 7.6]. In the critical case $p = (n+4+2\sigma)/(n-4)$, (4.1) is also guaranteed from [**JL14**, Theorem 3]. We can also select any regular radial solution to $(1.1)_\sigma$ in $\mathbf{R}^n$. Such a regular radial solution also enjoys (4.1), see [**FK19**]. In the supercritical case $p>(n+4+2\sigma)/(n-4)$, we cannot work directly with $U$. Instead, we are forced to make use of the Kelvin transform $\widehat{U}$ to go back to the subcritical case. In other words, in the supercritical case, instead of using $U$, we have to use $\widehat{U}$, which enjoys (4.1).

Thus, regardless of the size of $p$, we can always choose some $U$ satisfying the growth (4.1) and the following sharp pointwise inequality

$$-\Delta U \geqslant \alpha_{\rm GHWW}\frac{|\nabla U|^2}{U} \quad \text{in } \mathbf{R}^n\backslash\{0\}. \tag{4.2}$$

We shall use this in the proof of Lemma 4.2. Since this is straightforward, we omit the details.

Now we go back to the proof of the sharpness of the constant $2/(n-4)$ in Theorem 1.3. Our key point is that $U$ decays fast near infinity. In our previous work [NN26b], with help of the regularity of solutions near the origin, we can use the Kelvin transform to achieve the goal. In the underlying scenario, we need to choose the solution $U$ described above.

Since we already have a fast decay solution $U$ in hand, we can directly consider the Emden–Fowler transformation

$$V(t) = e^{(n-4)t} U(e^t) \quad \text{in } \mathbf{R}. \tag{4.3}$$

Our first task is to show that $V$ has the limit as $t \to +\infty$. In view of (4.1) we see that $V$ is bounded. Now making use of the inequality

$$(n-4)U(r) + rU'(r) \geqslant 0 \quad \text{in } (0, +\infty)$$

proved in Lemma 2.10 we know that $V$ is monotone. Hence, the following limit

$$V_\infty := \lim_{t \nearrow +\infty} V(t) \tag{4.4}$$

exists and is finite, where $r = e^t$. In addition, there holds $V_\infty > 0$, see Lemma 2.8. For convenience we denote $\tau = n-4$. Then

$$|\nabla U(x)| = |U'(r)| = \big|V'(t) - \tau V(t)\big| r^{-(\tau+1)}$$

and

$$\Delta U(x) = \big(V''(t) + (n-2\tau-2)V'(t) - \tau(n-\tau-2)V(t)\big) r^{-(\tau+2)}.$$

From these two identities we obtain from (4.2) the following

$$\begin{aligned} -V(t)V''(t) - (n-2\tau-2)V(t)V'(t) + \tau(n-\tau-2)V(t)^2 \\ \geqslant \alpha_{\rm GHWW}\big(V'(t) - \tau V(t)\big)^2 \end{aligned} \tag{4.5}$$

after canceling out the term $r^{-2\tau-2}$ from both sides. Then, we rewrite (4.5) in the following way

$$\begin{aligned} \big([\tau(n-\tau-2) - \alpha\tau^2]V(t) - (n-2\tau-2-2\alpha_{\rm GHWW}\tau)V'(t) - V''(t)\big)V(t) \\ \geqslant \alpha_{\rm GHWW}(V'(t))^2 \geqslant 0. \end{aligned}$$

We assume for a moment that

$$V^{(k)}(t) \to 0 \quad \text{as } t \to +\infty$$

for $k = 1, 2$. Then, together with (4.4), given arbitrary $\varepsilon > 0$ we choose large $T_\varepsilon \gg 1$ in such a way that

$$|(n-2\tau-2-2\alpha_{\rm GHWW}\tau)V'(t)| + |V''(t)| \leqslant \varepsilon \quad \text{for all } t \geqslant T_\varepsilon.$$

Then

$$\big([\tau(n-\tau-2) - \alpha_{\rm GHWW}\tau^2]V(t) + \varepsilon\big)V(t) \geqslant 0 \quad \text{for all } t \geqslant T_\varepsilon,$$

which then yields

$$[\tau(n-\tau-2) - \alpha_{\rm GHWW}\tau^2]V(t) + \varepsilon \geqslant 0 \quad \text{for all } t \geqslant T_\varepsilon.$$

From this together with (4.4) we should have

$$\tau(n-\tau-2) - \alpha_{\rm GHWW}\tau^2 \geqslant -\frac{\varepsilon}{2V_\infty} \quad \text{for large } t.$$

Since $\varepsilon > 0$ is arbitrary, we must have

$$\tau(n-\tau-2) \geqslant \alpha_{\rm GHWW}\tau^2,$$

which gives

$$\alpha_{\rm GHWW} \leqslant \frac{n-\tau-2}{\tau} = \frac{2}{n-4}.$$

Thus, we must have an equality, namely $\alpha_{\rm GHWW} = 2/(n-4)$ as claimed.

In the rest of the proof for this case, we show that $V^{(k)}(t) \to 0$ as $t \to +\infty$ with $k = 1, 2$. For clarity, we state and prove as lemmas. The reader is advised to refer to [**NN26b**] for further information.

Let us start with the limit of $V'$.

**Lemma 4.1.** *There holds $V'(t) \to 0$ as $t \to +\infty$.*

*Proof.* Since $V(t)$ has a limit as $t \to +\infty$, the conclusion follows provided $V''$ is bounded in $[1, +\infty)$. Easy calculation shows

$$V'(t) = \tau V(t) + r^{\tau+1} U'(r)$$

with $\tau = n-4$ as before. Making use of (2.15) we obtain

$$(n-2) r^\tau U(r) \geqslant -r^{\tau+1} U'(r) = r^{\tau+1} |U'(r)| \quad \text{for } r > 0.$$

Hence, with help of (4.1) we deduce that $r^{\tau+1}U'(r)$ is bounded, so is $V'$ because $V$ is bounded. Next, we have

$$V''(t) = (2\tau+1)V'(t) - \tau(\tau+1)V(t) + r^{\tau+2}U''(r).$$

Hence, it suffices to bound $r^{\tau+2}U''(r)$, and this can be proved by first making use of

$$r^{\tau+2}|U''(r)| \leqslant (n-1) r^{\tau+1} |U'(r)| + r^{\tau+2} |\Delta U(r)| \quad \text{for } r > 0.$$

Keep in mind that $U$ is monotone decreasing and $\Delta U$ is monotone increasing. Then, by integrating both sides of $(r^{n-1}U'(r))' = r^{n-1}\Delta U$ over $[1, r]$ we obtain

$$r^{n-1}U'(r) \leqslant U'(1) + \int_1^r s^{n-1} \Delta U(s) ds \leqslant U'(1) + \frac{r^n - 1}{n} \Delta U(r) \quad \text{for } r \geqslant 1.$$

Thanks to Lemma 2.7 we know that $\Delta U(r) = o(1)$ as $r \to +\infty$. Hence, there holds

$$U'(1) \leqslant \frac{\Delta U(r)}{n} < 0$$

for large $r$. Putting the preceding estimates together gives

$$r^{n-1}U'(r) \leqslant \frac{r^n}{n} \Delta U(r) \quad \text{for } r \gg 1.$$

Hence, again by (2.15) we eventually have

$$r^\tau U(r) \geqslant -\frac{1}{n-2} r^{\tau+1} U'(r) \geqslant -\frac{1}{n(n-2)} r^{\tau+2} \Delta U(r) > 0$$

for $r \gg 1$. This and (4.1) are enough for us to bound $r^{\tau+2}\Delta U(r)$, hence completing the present proof. □

Now we study the limit of $V''$.

**Lemma 4.2.** *There holds $V''(t) \to 0$ as $t \to +\infty$.*

*Proof.* Thanks to Lemma 4.1 and its proof we know that $V'$ and $V''$ is bounded in $[1, +\infty)$. Hence, to conclude the present proof, it suffices to prove that $V^{(3)}$ is bounded in $[1, +\infty)$. Without any difficulty, we quickly have

$$V^{(3)}(t) = 3(\tau+1)V''(t) - (3\tau^2 + 6\tau + 2)V'(t) + (\tau^3 + 3\tau^2 + 2\tau)V(t) + r^{\tau+3}U^{(3)}(r).$$

Hence, to bound $V^{(3)}$, it suffices to bound $r^{\tau+3}U^{(3)}(r)$. Notice that

$$(\Delta U)'(r) = U^{(3)}(r) + (n-1)r^{-1}U''(r) - (n-1)r^{-2}U'(r) \quad \text{for } r > 0,$$

which yields

$$0 \leqslant r^{\tau+3}|U^{(3)}(r)| \leqslant (n-1)r^{\tau+2}|U''(r)| + (n-1)r^{\tau+1}|U'(r)| + r^{\tau+3}|(\Delta U)'(r)|.$$

Therefore, it suffices to bound $r^{\tau+3}(\Delta U)'$, but this can be achieved via the following

$$r^{\tau+3}|(-\Delta U)'(r)| = r^{\tau+3}(\Delta U)'(r) \leqslant -(n-2)r^{\tau+2}\Delta U(r) \quad \text{for } r > 0,$$

thanks to (2.16), $\Delta^2 U > 0$, and the boundedness of $r^{\tau+2}\Delta U(r)$ established in the proof of the preceding lemma. This completes the present proof. □

The proof of Theorem 1.3 is complete.

## Ethical statements

**On the use of AI**. After announcing this work on arXiv:2610.39134, the authors used AI to audit the proof including checking calculation, logical issues, consistency, etc. After using this, the authors reviewed and edited the content as needed and take full responsibility for the content of the publication. In particular, this version does not contain AI written text.

**Data availability**. Data sharing is not applicable to this article as no dataset was generated or analysed during the current study.

## Acknowledgments

This research is funded by VNU University of Science under project number TN.26.24.

(Q.A. Ngô) University of Science, Vietnam National University, Hanoi, Vietnam, ORCID iD: 0000-0002-3550-9689

*Email address*: nqanh@vnu.edu.vn

(N.T. Trung) University of Science, Vietnam National University, Hanoi, Vietnam

*Email address*: trungisp58@gmail.com